\documentclass[12pt]{amsart}
\usepackage[cmtip,matrix,arrow,curve]{xy}
\usepackage{graphicx}
\def \fg{{\mathfrak g}}

\def \fk{{\mathfrak k}}

\def \dd{{\partial}}
\def \ii{\ensuremath{\mathbf{i}}}

\def \fa{{\mathfrak a}}

\def \fn{{\mathfrak n}}
\def \fq{{\mathfrak q}}

\def \su{{\mathfrak s}{\mathfrak u}}
\def \sl{{\mathfrak s}{\mathfrak l}}

\def \C{{\mathbb C}}

\def \R{{\mathbb R}}

\newtheorem{theorem}{Theorem}[section]
\newtheorem{lemma}[theorem]{Lemma}
\newtheorem{proposition}[theorem]{Proposition}

\numberwithin{equation}{section}

\begin{document}

\baselineskip=16pt

\title[The Poisson Geometry of SU(1,1)]{The Poisson Geometry of SU(1,1)}

\author[P. Foth]{Philip Foth}
\author[M. Lamb]{McKenzie Lamb}

\address{Department of Mathematics, University of Arizona, Tucson, AZ 85721}
\email{foth@math.arizona.edu}

\address{Mathematics and Computer Science, Ripon College, Ripon, WI 54971}
\email{lambm.ripon@gmail.com}



\date{August 14, 2009}

\begin{abstract}
We study the natural Poisson structure on the group ${\rm SU}(1,1)$ and related questions.
In particular, we give an explicit description of the Ginzburg-Weinstein isomorphism for
the sets of admissible elements. We also establish an analogue of Thompson's conjecture
for this group.
\end{abstract}

\maketitle
\section{Introduction}

The group ${\rm SU}(1,1)$ is the group of complex linear transormations of $\C^2$, which
preserves the pseudo-Hermitian form:
$$
\langle {\bf z},{\bf w}\rangle = z_1\bar{w}_1 - z_2\bar{w}_2\ .
$$
It is a subgroup of ${\rm SL}(2, \C)$, transversal to $AN$, consisting of upper-triangular
matrices with positive real diagonal entries. In the context of Poisson Lie groups, these
two Lie groups can be naturally viewed as dual to each other, and thus many questions
arise with regards to the induced Poisson structures.

In the present paper, we give explicit formulas in coordinates for those Poisson structures,
from which one can see the symplectic leaves. We also describe the dressing action and certain
natural identifications between subspaces of admissible elements in $AN$ and $\su(1,1)^*$.
In short, an element is called admissible if it has real spectrum, and its eigenvalue
corresponding to the timelike part is greater than the one for the spacelike part of $\C^2$.
This definition generalizes to all quasi-Hermitian Lie groups \cite{F071}, and beyond \cite{Neeb}.

In Sections 4 and 5 of the paper, we give two explicit approaches to a Poisson
isomorphism between the set of admissible elements in $AN$ with their
natural quadratic Poisson structure and the corresponding set in $\su(1,1)^*$ with the Lie-Poisson structure.
This can be considered as the first step in generalizing the Ginzburg-Weinstein theorem \cite{GW}
to the non-compact setup. One of our approaches follows the original path by Ginzburg and Weinstein,
and the other follows the idea of Flaschka and Ratiu \cite{FR}, based on the Gelfand-Tsetlin
coordinates. In fact, their conjecture was proven by Alekseev and Meinrenken \cite{AM} for
the general ${\rm SU}(n)$ case. However, in all of those approaches, compactness was used
quite heavily, so that a direct generalization to a pseudo-unitary setup is not possible. This is
one of the reasons that we need to be explicit in our constructions.

In the last section, we establish an analogue of Thompson's conjecture in the pseudo-unitary setup associated to the group ${\rm SU}(1,1)$.  Thompson's conjecture concerns the equality of spectra in the linear and non-linear situations, for the sums and products of admissible elements, respectively. In our
case, the singular spectrum is replaced by the so-called \emph{admissible spectrum}, related to the
decomposition of an open subset of admissible elements in a complex reductive group $G_\C$ as the product
$HA_{\rm adm}H$, where $H$ is a quasi-Hermitian real form of $G_\C$.

\section{Basic Facts}

The group $G={\rm SU}(1,1)$ is realized as the group of $2\times 2$ complex matrices
$$
\left( \begin{array}{cc} u & v \\ \bar{v} & \bar{u} \end{array} \right)\
$$
of determinant equal to $1$, i.e. satisfying $|u|^2 - |v|^2 = 1$. This group is
isomorphic to the group ${\rm SL}(2, \R)$ as a real Lie group. The Lie algebra $\fg$
of $G$ is
$$
\fg = \left\{ \left( \begin{array}{cc}
\ii r & \eta \\ \bar{\eta} & -\ii r \end{array} \right) : \ r\in\R, \ \eta\in\C \right\}.
$$

A convenient basis for $\fg$ consists of the elements $X$, $Y$, and $H$, where
$$
X= \left( \begin{array}{cc} 0 & 1 \\ 1 & 0 \end{array} \right), \
Y= \left( \begin{array}{cc} 0 & \ii \\ -\ii & 0 \end{array} \right), \ \ {\rm and} \ \
H= \left( \begin{array}{cc} \ii & 0 \\ 0 & -\ii \end{array} \right).
$$
The Lie algebra structure on $\fg$ as well as the Poisson structure on $\fg^*$ is generated
by the Lie bracket relations:
$$
[X,Y]=2H, \ \ [X,H] = - 2Y, \ \ [Y,H] = 2X\ .
$$
The dual vector space $\fg^*$ can be identified with the subspace of $\sl(2, \C)$ of the form
$$
\fg^* = \left\{ \left( \begin{array}{cc}
z & x+\ii y \\ -x+\ii y & -z \end{array} \right) : \ x,y,z\in\R \right\},
$$
consisting of pseudo-Hermitian matrices of signature $(1,1)$. The natural pairing between
$\fg$ and $\fg^*$ is given by the non-degenerate form
$$
\langle A, B \rangle = {\Im}({\rm Tr}(AB))\ .
$$

The linear Poisson structure on $\fg^*$ is then given by
\begin{equation}
\pi_{0} = -z \frac{\partial}{\partial x}\wedge\frac{\partial}{\partial y} +
y \frac{\partial}{\partial z}\wedge\frac{\partial}{\partial x} +
x \frac{\partial}{\partial y}\wedge\frac{\partial}{\partial z} \ .
\label{eq:epi0}
\end{equation}
A Casimir for $\pi_0$ is $z^2 - x^2 - y^2$.

Let $\fa$ and $\fn$ be the Lie subalgebras of $\fg_\C = \sl(2, \C)$ of the form
$$
\fa = \left\{ \left( \begin{array}{cc}
 r & 0 \\ 0 & - r \end{array} \right) : \ r\in\R \right\}, \ \ \
\fn = \left\{ \left( \begin{array}{cc}
 0 & n \\ 0 & 0 \end{array} \right) : \ n\in\C \right\}.
$$

The subalgebra $\fa+\fn$ of $\fg_\C$ can also be viewed as the dual vector space for
$\fg$ and the natural pairing is given by
$$
\langle A, B \rangle = 2{\Im}({\rm Tr}(AB))\ .
$$

Thus, we have a Manin triple, i.e., a pair of transversal, Lagrangian subalgebras, $\fg$ and $\fa+\fn$
of $\sl(2, \C)$, together with the above non-degenerate pairing. We refer to \cite{LW} for preliminaries on
Poisson Lie groups. We therefore have induced Poisson Lie group
structures $\pi_G$ and $\pi_{AN}$ on $G$ and $AN= \exp(\fa+\fn)$ respectively.

The Poisson structure $\pi_G$ can be expressed as follows. Consider the element
$\Lambda\in\fg\wedge\fg$ given by $\Lambda = \frac{1}{2}X\wedge Y$. Then for any $g\in G$ we have
$$
\pi_G(g) = (r_g)_*\Lambda - (l_g)_*\Lambda\ ,
$$
where $r_g$ and $l_g$ are the right and left translations respectively by $g$ on $G$.

In terms of the matrix elements, for
$g=\displaystyle{\left( \begin{array}{cc} u & v \\ \bar{v} & \bar{u} \end{array} \right)}$,
this bracket is quadratic and given by:
$$
\{ u, \bar{u}\} = -2\ii |v|^2 \ \ \ \ \ \ \ \
\{ u, v\} = -\ii uv \ \ \ \ \ \ \ \
\{ u, \bar{v}\} = -\ii u\bar{v} \ \ \ \ \ \ \ \
$$
$$
\{ \bar{u}, v\} = \ii \bar{u}v \ \ \ \ \ \ \ \
\{ \bar{u}, \bar{v}\} =  \ii \bar{u}\bar{v} \ \ \ \ \ \ \ \
\{ v, \bar{v}\} = 0 \ \ \ \ \ \ \ \
$$

The matrix formula for the Poisson structure $\pi_{AN}$ on
$$
AN = \left\{ \left( \begin{array}{cc}
 \exp(\frac{z}{2}) & x+\ii y \\ 0 & \exp(-\frac{z}{2}) \end{array} \right) : \ x,y,z\in\R \right\}
 $$
is as follows:
$$
\pi_{AN} = - \sinh(z) \frac{\partial}{\partial x}\wedge\frac{\partial}{\partial y} +
y \frac{\partial}{\partial z}\wedge\frac{\partial}{\partial x} +
x \frac{\partial}{\partial y}\wedge\frac{\partial}{\partial z}\ .
$$
A Casimir function for this Poisson structure is
$$
\xi(x,y,z) = 2\cosh (z) - x^2 - y^2\ .
$$

The linearization of $\pi_{AN}$ is denoted by $\pi_0$.  This should not lead to
confusion because it is given by the exact same formula, (\ref{eq:epi0}).

Consider the complex anti-linear involution $\dagger$ on the Lie algebra
$\fg_\C = \sl(2, \C)$, given by
$$
M^\dagger = J\bar{M}^TJ, \ \ {\rm where} \ \
J = \left( \begin{array}{cc} 1 & 0 \\ 0 & -1 \end{array} \right) .
$$
Let $\fq$ denote the fixed point subspace of this involution. (Note that
$\fg=\su(1,1)$ is the $(-1)$-eigenspace of $\dagger$.) Earlier, we have identified
$\fq$ with $\fg^*$. The same formula defines an involution on the matrix Lie group
${\rm SL}(2, \C)$, and its fixed point set there is denoted by
$$
Q = \left\{ \left( \begin{array}{cc} c & \beta \\ - \bar{\beta} & d \end{array}
\right), \ \ c,d\in\R, \ \ \beta\in\C, \ \ cd+|\beta|^2 = 1 \right\} \ .
$$

\section{Admissible loci}

For an open subset $Q'$ of $Q$, on which $c\ne 0$, one can define coordinates:
$$
Q' = \left\{ \left( \begin{array}{cc} c & a+\ii b \\ - a+\ii b &
\frac{1-a^2-b^2}{c} \end{array}
\right), \ \ a,b,c,d\in\R \right\} \ .
$$
Consider the following symmetrization map:
$$
{\rm Sym}: AN\to Q, \ \ M\mapsto M^{\dagger}M\ , \ \ {\rm or}
$$
$$
a = x e^{z/2}, \ \ \ b = y e^{z/2}, \ \
c = e^z
$$
in coordinates.
Under this map, the Poisson tensor $\pi_{AN}$ pushes down to
$$
\pi_Q = \frac{1}{2}(1-a^2-b^2-c^2)\frac{\partial}{\partial a}\wedge\frac{\partial}{\partial b}+
bc\frac{\partial}{\partial c}\wedge\frac{\partial}{\partial a}+
ac\frac{\partial}{\partial b}\wedge\frac{\partial}{\partial c}\ ,
$$
with Casimir
$$
F(a,b,c) = \frac{1+c^2-a^2-b^2}{c}\ ,
$$
which is simply the trace.  Note that the dressing action of $G$ on $AN$ converts to conjugation on $Q$.

We define the subset of admissible elements $\fq_{\rm adm}\subset\fq$ as the set of elements
conjugate to the diagonal matrices ${\rm diag}(\lambda, -\lambda)$ with $\lambda > 0$.
The set of adimissible elements forms an open cone in $\fq$ defined by
$$
z^2-x^2-y^2 > 0, \ \ {\rm and} \ \ z > 0.
$$
Denote also $Q_{\rm adm} = \exp(\fq_{\rm adm})$. The exponential map is easily checked to
be invertible and we denote its inverse
$$
\log: \ \ Q_{\rm adm}\to \fq_{\rm adm} \ .
$$
Actually, if one denotes $A_{\rm adm} = {\rm diag}(e^{z/2}, e^{-z/2})$, and
$(AN)_{\rm adm}=G.A_{\rm adm}$, where ``." denotes the dressing action,
then the image of $(AN)_{\rm adm}$ under the symmetrization map is exactly
$Q_{\rm adm}$. On the set of admissible elements, $(AN)_{\rm adm}$, the
right dressing action is globally defined. (Which is not true for the left
dressing action.) This is also the case for a general pseudo-unitary groups
${\rm SU}(p,q)$. In our case, we have explicitly:
$$
\left( \begin{array}{cc} e^{z/2} & 0 \\ 0 & e^{-z/2} \end{array} \right)\cdot
\left( \begin{array}{cc} u & v \\ \bar{v} & \bar{u} \end{array} \right) =
\left( \begin{array}{cc} u' & v' \\ \bar{v}' & \bar{u}' \end{array} \right) \cdot
\left( \begin{array}{cc} \rho & m \\ 0 & \rho^{-1} \end{array} \right) \ ,
$$
where
$$
\rho = \sqrt{|u|^2e^z - |v|^2e^{-z}}, \ u' = \frac{ue^{z/2}}{\rho}, \
v' = \frac{ve^{-z/2}}{\rho}, \ m = 2\bar{u}v\frac{\sinh(z)}{\rho}\ .
$$
Note that $\rho$ is well-defined, since $z>0$.
It is also easy to see that the symmetrization map is a diffeomorphism on the
set of admissible elements.

For an element $B\in Q_{\rm adm}$, denote by $(e^{\lambda}, e^{-\lambda})$ the set of its
eigenvalues with $\lambda >0$. Then the log map from $Q_{\rm adm}$ to $\fq_{\rm adm}$,
$$
{\rm log}: \ \ \left( \begin{array}{cc} c & a+\ii b \\ - a+\ii b &
\frac{1-a^2-b^2}{c} \end{array}
\right) \ \mapsto \ \left( \begin{array}{cc}
z & x+\ii y \\ -x+\ii y & -z \end{array} \right)
$$
is given by
$$
(a,b,c) \ \mapsto \
\left( a \frac{\lambda}{\sinh(\lambda)}, b \frac{\lambda}{\sinh(\lambda)},
c \frac{\lambda}{\sinh(\lambda)} - \lambda\coth(\lambda)\right) \ = \ (x,y,z)\ .
$$
Under the composition of these two diffeomorphisms on the sets of admissible elements,
the Poisson structure $\pi_{AN}$ pushes down to
\begin{eqnarray*}
\pi &:=& \log_*({\rm Sym}_*(\pi_{AN}))\\
&=& - z (\lambda\coth\lambda+z) \frac{\partial}{\partial x}\wedge\frac{\partial}{\partial y} +
y (\lambda\coth\lambda+z) \frac{\partial}{\partial z}\wedge\frac{\partial}{\partial x} +
x (\lambda\coth\lambda+z) \frac{\partial}{\partial y}\wedge\frac{\partial}{\partial z}\\
&=& (\lambda\coth\lambda+z)\pi_0\ .
\end{eqnarray*}
The natural Casimir function for this Poisson structure is clearly
$x^2+y^2-z^2 = -\lambda^2$, which is the determinant. Note that the symplectic leaves of those two
structures are actually the same, and are hyperboloids. One can view $\pi$ as a family of Poisson
structures on a single hyperboloid, depending on $\lambda$, and identify it diffeomorphicaly with the
lower hemisphere, and show that it extends to the whole sphere. In terms of a holomorphic coordinate
$w$ on the sphere, this one-parameter family of Poisson structures can also be written as
$$
\pi(\tau) = \ii(1-|w|^2)|w|^2\ \frac{\dd }{\dd w}\wedge\frac{\dd }{\dd \bar{w}} +
\tau\cdot \ii(1-|w|^2)^2\ \frac{\dd }{\dd w}\wedge\frac{\dd }{\dd \bar{w}} ,
$$
for $\tau\in\R$, where the first term is the so-called $\Pi_v$ structure from
\cite{FL2}, and the second is
an ${\rm SU}(1,1)$-invariant Poisson structure on $S^2$.

\section{Poisson isomorphism: the Ginzburg-Weinstein approach}

The Ginzburg-Weinstein approach in the compact situation to finding a Poisson isomorphism between
$\fk^*$ and $K^*$ was to prove the existence of a vector field whose flow would connect
$\pi_0$ and $\pi$.
We will construct such a vector field explicitly for the $SU(1,1)$ case.

Following the Ginzburg-Weinstein argument, define a bivector field $\pi_t$ on
$\fq$ by $$\pi_t(\vec{v}):=\frac{\pi(t\vec{v})}{t},$$
where the expression on the right-hand side is identified with an element of
$\wedge^2\left(T_{\vec{v}}\fq\right)$ by translation.  Now set $$\dot{\pi}_t:=\frac{d}{dt}\pi_t$$ and
$$\dot{\pi}:=\left.\frac{d}{dt}\right|_{t=1}\pi_t.$$
In coordinates,
\begin{eqnarray*}\dot{\pi}
&=&-z\left(\lambda\coth\lambda+z-\frac{\lambda^2}{\sinh^2\lambda}\right)
\frac{\dd}{\dd{x}}\wedge\frac{\dd}{\dd{y}} +\\
&&y\left(\lambda\coth\lambda+z-\frac{\lambda^2}{\sinh^2\lambda}\right)
\frac{\dd}{\dd{z}}\wedge\frac{\dd}{\dd{x}} + \\
&&x\left(\lambda\coth\lambda+z-\frac{\lambda^2}{\sinh^2\lambda}\right)
\frac{\dd}{\dd{y}}\wedge\frac{\dd}{\dd{z}}.
\end{eqnarray*}
\begin{proposition}
\label{ExistenceOfYProp} There exists a vector field $X$ on $\fq_{adm}\cong(\fa+\fn)_{adm}$ such that
\begin{itemize}
\item[(1)] $[X,\pi]=\dot{\pi}$,
\item[(2)] $X$ has the zero linearization at the origin,
\item[(3)] \label{YTangent}$X$ is tangent to the symplectic leaves of $\pi$ (and $\pi_0$), and
\item[(4)] $X$ is complete.
\end{itemize}
\end{proposition}
Note that in the compact case considered by Ginzburg and Weinstein, completeness
simply follows from the fact that $X$ is tangent to the symplectic leaves of $\pi$.

\begin{proof}
It will be convenient to convert to \emph{hyperbolic coordinates} $(\lambda,\phi,s)$.
The relations between rectangular and hyperbolic coordinates are
$$\begin{array}{rclcrcl}
x&=&\lambda(\sinh s)(\cos\phi) &\hskip.5in & \lambda&=&\sqrt{z^2-x^2-y^2}\\
y&=&\lambda(\sinh s)(\sin\phi) &\hskip.5in & \phi&=&\arctan (y/x)\\
z&=&\lambda(\cosh s) &\hskip.5in & s&=&\cosh^{-1}(z/\lambda).
\end{array}$$

In these coordinates,
$$\pi=\frac{1}{\sinh s}\left(\coth \lambda + \cosh s\right)
\frac{\dd}{\dd{\phi}}\wedge\frac{\dd}{\dd{s}},$$ and
$$\dot{\pi}=\frac{1}{\sinh s}\left(\coth \lambda + \cosh s-\frac{\lambda}{\sinh^2\lambda}\right)
\frac{\dd}{\dd{\phi}}\wedge\frac{\dd}{\dd{s}}.$$

Set
\begin{eqnarray*}g(s)&:=&\frac{1}{\sinh s}\left(\coth \lambda + \cosh s\right),\\
h(s)&:=&\frac{1}{\sinh s}\left(\coth \lambda + \cosh s-\frac{\lambda}{\sinh^2\lambda}\right).
\end{eqnarray*}
The action of the diagonal torus $T\subset G$ on $\fa+\fn$ corresponds to rotation about
the $z$-axis in $\fq$.  Since $\pi$, and hence $\dot{\pi}$, are invariant under the
torus action, we may assume that $X$ has the form $\displaystyle{f(s)\frac{\dd}{\dd{s}}}$,
where $f(s)$ does not depend on $\phi$.  Then the equation $[X,\pi]=\dot{\pi}$ reduces to the ODE
$$f\frac{\partial g}{\partial s}-g\frac{\partial f}{\partial s}=h.$$
Rewriting the left-hand side using the quotient rule gives
$\displaystyle{ -\frac{\dd}{\dd{s}}\left(\frac{f}{g}\right)\cdot g^2=h}$, or, equivalently,
\begin{equation}\label{ODEIntegral}f=-g\cdot\int\frac{h}{g^2}\;ds\ .\end{equation}
Integrating, we obtain
\begin{equation}\label{YIntegral}\int\frac{h}{g^2}\;
ds=\ln(\coth \lambda+\cosh s)+\frac{\lambda}{\sinh^2\lambda}\left(\frac{1}{\coth \lambda+\cosh s}\right)+C,\end{equation}
where $C$ is constant with respect to $s$ and $\phi$.
Note that $g\to\infty$ as $s\to 0$.  Therefore, to ensure smoothness when $s=0$, set
$$
C= - \ln(\coth
\lambda+1)-\frac{\lambda}{\sinh^2\lambda}\left(\frac{1}{\coth \lambda+1}\right)\ .
$$
Thus, we obtain the vector field
$$
X=
$$
$$
-\left(\frac{\coth \lambda + \cosh s}{\sinh s}\right)\cdot\Bigg[\ln\left(\frac{\coth
\lambda+\cosh s}{\coth \lambda+1}\right)+\frac{\lambda}{\sinh^2\lambda}\left(\frac{1}{\coth
\lambda+\cosh s}-\frac{1}{\coth \lambda+1}\right)\Bigg]\frac{\dd}{\dd{s}},
$$
which extends smoothly to the positive $z$-axis (where it vanishes).
This vector field is smooth on the open cone $z>\sqrt{x^2+y^2}$, extends
continuously to the boundary $z=\sqrt{x^2+y^2}$, and satisfies $[X,\pi]=\dot{\pi}$.

It is easy to check that $X$ has zero linearization at the origin. Since when $s\to\infty$, we have
$\displaystyle{X\sim s\frac{\dd}{\dd{s}}}$,
 the restriction $X_{\lambda}$ of $X$ to any hyperboloid $\lambda=\sqrt{z^2-x^2-y^2}$
extends continuously to the boundary, which we have identified with the unit circle
in the plane.  Since the closed unit disk is compact, it follows that $X_{\lambda}$
is complete for every $\lambda$, which implies that $X$ is complete.  This completes the proof.
\end{proof}

Given the vector field $X$ from Proposition \ref{ExistenceOfYProp}, the
Ginzburg-Weinstein argument goes through as in the compact case.  Defining
$X_t$ by $$X_t(\vec{p}):=\frac{X(t\vec{p})}{t^2},$$
the corresponding flow $\varphi_t$ pushes $\pi_0$ forward to $\pi_t$, and
in particular, $\left(\varphi_1\right)_*\left(\pi_0\right)=\pi$.
Thus, $\phi_1$ is the desired Poisson isomorphism.

\section{Poisson Isomorphism: The Flaschka-Ratiu Approach}
In this section, we apply the procedure used by Flaschka and Ratiu in \cite{FR}
for the ${\rm SU}(2)$ case to our pseudo-unitary situation.
The idea is to use the Gelfand-Tsetlin coordinates, which are given by the eigenvalues
of the principal minors of a matrix.  The Gelfand-Tsetlin coordinates were developed in \cite{GS} for the unitary case and extended in \cite{F071} to the pseudo-unitary case.

The elements of $\fq_{\rm adm}$ with eigenvalues $(\lambda, -\lambda)$ for $\lambda >0$
can be parameterized by the matrices
\begin{equation}\label{FraqParam1}\begin{pmatrix}z & \sqrt{z^2-\lambda^2}\cdot e^{\ii\theta}\\
-\sqrt{z^2-\lambda^2}\cdot e^{-\ii\theta} & -z.\end{pmatrix},\end{equation}
with $z\ge \lambda$ and $0\leq\theta<2\pi$.
Define coordinates on $\fq$ by
identifying (\ref{FraqParam1}) with $(z,\lambda,\theta)\in\R^3$.
The coordinates $z$ and $\lambda$ are the eigenvalues of the upper left $1\times1$ and $2\times2$ minors.
In these coordinates, $$\pi_0=\frac{\dd}{\dd{\theta}}\wedge\frac{\dd}{\dd{z}}.$$
The symplectic structure induced by $\pi_0$ on \emph{any} symplectic
leaf
$$
\Theta_{\lambda}=\left\{ (x,y,z)\;:\; \sqrt{z^2-y^2-x^2}=\lambda\right\}
$$
is then given by:
$$\omega_0:=d\theta\wedge dz.$$
Similarly, the elements of $Q$ with eigenvalues $(e^\lambda, e^{-\lambda})$ can be parameterized by the matrices
\begin{equation}\label{FraqParam2}\begin{pmatrix}e^{w} & \sqrt{(e^{w}-e^{\lambda})(e^{w}-e^{-\lambda})}\cdot e^{\ii\theta}\\
-\sqrt{(e^{w}-e^{\lambda})(e^{w}-e^{-\lambda})}\cdot e^{-\ii\theta} & 2\cosh(\lambda)-e^{w}\end{pmatrix}.\end{equation}
Define coordinates on $Q$ by
identifying (\ref{FraqParam1}) with $(w,\lambda,\theta)\in\R^3$.
In these coordinates: $$\pi_Q=\frac{\dd}{\dd{\theta}}\wedge\frac{\dd}{\dd{w}}.$$
The symplectic structure induced by $\pi_Q$ on any symplectic leaf
$$
\Psi_{\lambda}=\left\{(a,b,c)\;:\;\frac{1+c^2-a^2-b^2}{c}=2(\cosh\lambda)\right\}
$$
is then given by:
$$\omega_{Q}:=d\theta\wedge dw.$$
Given these simple expressions for $\omega_0$ and $\omega_1$, for each $\lambda$, we can define a symplectomorphism
from $\Theta_{\lambda}$ to $\Psi_{\lambda}$ by identifying the matrix
(\ref{FraqParam1}) with the matrix (\ref{FraqParam2}).
A Poisson isomorphism $f$ from
$(\fq_{\rm adm},\pi_0)$ to $(Q_{\rm adm},\pi_Q)$ is
obtained by allowing $\lambda$ to vary over the interval $(0,\infty)$.
Equivalently, $f$ sends $(z,\lambda,\theta)$ to $(w,\lambda,\theta)$.  In terms of the coordinates
$$\begin{pmatrix}
z & x+iy\\
-x+iy & -z\end{pmatrix}\leftrightarrow(x,y,z)$$ on $\fq$ and
$$\left(\begin{array}{cc}
c & a+ib\\
-a+ib & \frac{1-a^2-b^2}{c}\end{array}\right)\leftrightarrow(a,b,c)$$ on $Q$, $f$ is given by
\begin{eqnarray*}
a&=&\sqrt{e^{2z}-2e^{z}\cosh(\sqrt{z^2-x^2-y^2})+1}\left(\frac{x}{\sqrt{x^2+y^2}}\right)\\
b&=&\sqrt{e^{2z}-2e^{z}\cosh(\sqrt{z^2-x^2-y^2})+1}\left(\frac{y}{\sqrt{x^2+y^2}}\right)\\
c&=&e^{z}.
\end{eqnarray*}
The map $f$ is neither one-to-one nor onto, but it is actually a diffeomorphism
when restricted to $\fq_{\rm adm}$.

\section{Thompson's conjecture}

In this section we establish an analogue of Thompson's conjecture
for the ${\rm SU}(1,1)$ case. The original proof
in the unitary case is due to Klyachko \cite{Kl}. First, we say that an element
$g\in G_\C = {\rm SL}(2, \C)$ is admissible if it can be decomposed as a product $g=hb$,
with $h\in G={\rm SU}(1,1)$ and $b\in(AN)_{\rm adm}$. In this case, the {\it admissible spectrum}
of $g$ is the number $\gamma > 0$ such that the pair $(e^{\gamma}, e^{-\gamma})$ is the spectrum
of $g^{\dagger}g=b^{\dagger}b$.

This is equivalent to saying that $g\in G_\C$ is admissible if it lies in the open subset
$GA_{\rm adm}G$ of $G_\C$, in which case its admissible spectrum can be read from the middle
term of this decomposition.

For example, if
\begin{equation}
b = \left( \begin{array}{cc} e^{z/2} & x+\ii y \\ 0 & e^{-z/2} \end{array}\right) \ \in \ (AN)_{\rm adm}\ ,
\label{eq:eb144}
\end{equation}
then it is admissible if and only if $z>0$ and $\Delta:=e^{z} + e^{-z} - (x^2 +y^2) >2$, in which case the
admissible spectrum is given by
$$
{\gamma} = \log \left( \frac{\Delta + \sqrt{\Delta^2 - 4}}{2}\right)\ .
$$

\begin{lemma}
If $g_1, g_2 \in G_{\C}$ are admissible, then their product $g_1g_2$ is also admissible.
\end{lemma}

\begin{proof} The proof is omitted, as it is a short computational affair,
which uses the fact that the dressing action on $(AN)_{\rm adm}$ does not change
the admissible spectrum, and therefore we can assume one of the two elements
in $(AN)_{\rm adm}$ is diagonal.
\end{proof}

It is easy to see that the possible admissible spectrum of an element $b$ given by (\ref{eq:eb144})
lies in the interval $[z, \infty)$.

Next, one can readily establish that for two elements $M_1, M_2\in \fq_{\rm adm}=\su(1,1)^*_{\rm adm}$
with respective eigenvalues $(\lambda_1, -\lambda_1)$ and $(\lambda_2, -\lambda_2)$ such that
$\lambda_1, \lambda_2 > 0$,
the possible spectrum $(\lambda, -\lambda)$ of $M_1+M_2$ satisfies $\lambda\ge \lambda_1+\lambda_2$,
which is equivalent to the reversed triangle inequality in Minkowski space \cite{F081}.

Thompson's conjecture in our particular case now is equivalent to the following
\begin{proposition}
For two admissible elements $g_1$ and $g_2$ from $G_\C$ with admissible spectra $\lambda_1$ and
$\lambda_2$ respectively, the admissible spectrum of their product $g_1g_2$ lies in the interval
\ $[\lambda_1+\lambda_2, \infty)$.
\end{proposition}
\begin{proof}

Clearly, we can assume $g_1=b_1\in(AN)_{\rm adm}$ and $g_2=b_2\in(AN)_{\rm adm}$ as well.
Also, using the dressing action, one can assume that one of those elements, say $b_2$, is diagonal:
$b_2=a_2\in A_{\rm adm}$, where \ $a_2 = {\rm diag}(\rho, \rho^{-1})$.
Let $a_1={\rm diag}(r,r^{-1})$ and $g\in {\rm SU}(1,1)$ be such that $b_1=a_1.g$, with respect
to the dressing action. The element $g$ is given by
$$
g = \left( \begin{array}{cc} u & v \\ \bar{u} & \bar{v} \end{array} \right),
$$
with $|u|^2-|v|^2=1$. The admissible spectra of $b_1$ and $b_2$ are
$\lambda_1=2\log(r)$ and $\lambda_2=2\log(\rho)$ respectively.
Now consider the product $b=b_1b_2$ and compute:
$$
b^{\dagger}b = b_2^{\dagger}b_1^{\dagger}b_1b_2 =
a_2^{\dagger} (a_1.g)^{\dagger}(a_1.g) a_2 =
a_2g^{-1}a_1^2ga_2\ .
$$
In terms of matrices, we have
$$
{\rm Tr}(b^{\dagger}b) = r^2 \rho^2|u|^2 - r^{-2}\rho^2|v|^2+r^{-2}\rho^{-2}|u|^2
-r^2\rho^{-2}|v|^2\ .
$$
Thus, if $\mu$ is the greatest root of the quadratic equation
$$
\mu+\frac{1}{\mu}={\rm Tr}(b^{\dagger}b)\ ,
$$
then the admissible spectrum of $\lambda$ of $b=b_1b_2$ is given by $\lambda=\log(\mu)$.
It follows that Thompson's conjecture in our case is equivalent to proving
$$
\mu\ge r^2\rho^2\ .
$$
Now, if we substitute $|u|^2 = 1+|v|^2$, then we obtain:
$$
\mu+\frac{1}{\mu} =
|v|^2\left( r^2\rho^2-\frac{\rho^2}{r^2}+\frac{1}{r^2\rho^2}-\frac{r^2}{\rho^2}\right)
+r^2\rho^2+\frac{1}{r^2\rho^2} =
$$
$$
= |v|^2 \left( \rho^2-\frac{1}{\rho^2}\right) \left( r^2-\frac{1}{r^2}\right)
+r^2\rho^2+\frac{1}{r^2\rho^2} \ge r^2\rho^2+\frac{1}{r^2\rho^2}\ .
$$
Since $\mu$ is the greater root of this equation, and $\displaystyle{\xi(x) = x+\frac{1}{x}}$
is an increasing function of $x$, for $x>1$, we conclude that
$\mu\ge r^2\rho^2$ as desired.
\end{proof}

%


\end{document}